\documentclass[preprint,12pt]{elsarticle}

\usepackage[english]{babel}
\usepackage[left=1in,right=1in,top=1in,bottom=1in]{geometry}
\usepackage{amsmath,amssymb,amsthm}
\usepackage{url}
\usepackage{graphicx}
\usepackage{float}
\newtheorem{theorem}{Theorem}[section]
\newtheorem{lemma}[theorem]{Lemma}
\newtheorem{proposition}[theorem]{Proposition}
\newtheorem{corollary}[theorem]{Corollary}
\theoremstyle{definition}

\newtheorem{definition}[theorem]{Definition}

\newtheorem{remark}[theorem]{Remark}
\numberwithin{equation}{section}

\journal{Journal of Mathematical Analysis and Applications}

\begin{document}
	
	\begin{frontmatter}
		
		
		
		\title{Stability Analysis of Grünwald Interpolation Operators on Chebyshev Nodes} 
		
		
		\author{P. C. Vinaya} 
		
		\affiliation{organization={Cochin University of Science and Technology},
			addressline={Kalamassery}, 
			city={Cochin},
			postcode={682022}, 
			state={Kerala},
			country={India}}
		
		\begin{abstract}
			In 1941, G.~Grünwald proved the convergence of a sequence of operators constructed using classical Lagrange interpolation at Chebyshev nodes. In this work, we establish a perturbed version of Grünwald’s result, thereby extending the class of admissible nodal points. Specifically, we provide sufficient conditions for convergence when the interpolation nodes are of the form $\{\cos \eta_k\}_{k=1}^n$, where $\{\eta_k\}_{k=1}^n$ is a general sequence. We refer to these operators as Grünwald operators. In particular, we prove a convergence result when $\{\eta_k\}_{k=1}^n$ is equidistant and uniformly distributed. We establish a Voronovskaja-type estimate for the convergence of these operators and derive quantitative results using modulus of continuity.
		\end{abstract}

		\begin{keyword}
			Lagrange interpolation, Chebyshev nodes, Voronovskaja-type theorem, rate of convergence, modulus of smoothness.
			\MSC[2020]  41A25\sep 41A60\sep 41A17\sep 41A35\sep 41A50
			
		\end{keyword}
		
	\end{frontmatter}
	
	\section{Introduction} 
	Given a continuous function \( f \in C[-1,1] \), the Lagrange interpolation operator provides a unique polynomial of degree \( n-1 \) that interpolates the values \( f(x_1), f(x_2), \ldots, f(x_n) \) at the nodes \( x_1, x_2, \ldots, x_n \), respectively. It is expressed as
	\[
	\mathcal{L}_n(f)(x) = \sum_{k=1}^{n} f(x_k) P_k(x), \quad x \in [-1,1],
	\]
	where the Lagrange basis (or fundamental) polynomials \( P_k(x) \) are defined by
	\[
	P_k(x) = \frac{\zeta(x)}{\zeta'(x_k)(x - x_k)}, \quad k = 1, \ldots, n,
	\]
	with
	\[
	\zeta(x) = c(x - x_1)(x - x_2)\cdots(x - x_n),
	\]
	and \( c \) is an arbitrary non-zero constant. It is easy to verify that
	$
	\sum\limits_{k=1}^n P_k(x) = 1$ for all $x \in [-1,1]$. The Chebyshev nodes of the first kind are given by
	\[
	x_k = \cos \theta_k^{(n)}, \quad \text{where} \quad \theta_k^{(n)} = \frac{(2k - 1)\pi}{2n}, \quad k = 1, 2, \dots, n.
	\]
	
	Using this, the Lagrange interpolation operator based on Chebyshev nodes can be expressed as
	\[
	\mathcal{L}_n(f)(\theta) = \sum_{k=1}^n f(\cos \theta_k^{(n)}) \cdot (-1)^{k+1} \frac{\cos(n\theta) \sin \theta_k^{(n)}}{n(\cos\theta - \cos \theta_k^{(n)})},
	\]
	where
	\[
	P_k(\theta) = (-1)^{k+1} \frac{\cos(n\theta) \sin \theta_k^{(n)}}{n(\cos\theta - \cos \theta_k^{(n)})},
	\quad k = 1, \dots, n.
	\] 
	The convergence of the Lagrange interpolation operators is a long-standing problem in the literature~\cite{erdos, erdos1, nevai}. 
	It is well known that the Lagrange interpolation operators \( \mathcal{L}_n(f) \) do not converge to \( f \) at all points in the interval \( [-1,1] \)~\cite{faber}.

	In 1936 and 1937, G.~Grünwald \cite{grunwald1} and Józef Marcinkiewicz \cite{josef} independently addressed the problem of convergence of Lagrange interpolation operators on Chebyshev nodes. The celebrated Grünwald–Marcinkiewicz theorem asserts the following:
	
	\begin{theorem}[Grünwald–Marcinkiewicz] \cite{grunwald1,josef}
		There exists a function $f\in C[-1,1]$ such that for all $x\in [-1,1]$, the sequence $\mathcal{L}_n(f)(x)$ diverges.	
	\end{theorem}
	A detailed exposition of the Grünwald–Marcinkiewicz theorem, along with insights into the lives of these remarkable mathematicians, can be found in \cite{mills, vertesi}.\par
	We define a sequence of operators using the Lagrange interpolation operators below. 
	\begin{definition}\label{Gn}
		Let \( \mathcal{G}_n: C[0,\pi] \to C[0,\pi] \) be an operator defined by
		\[
		\mathcal{G}_n(f)(\theta) = \frac{1}{2} \sum_{k=1}^n f(\eta_k) \left\{ l_k\left(\theta - \frac{\pi}{2n}\right) + l_k\left(\theta + \frac{\pi}{2n}\right) \right\},
		\]
		where \( l_k \) are the fundamental polynomials of the Lagrange interpolation operator defined on the nodes \( \{ \cos \eta_k \}_{k=1}^n \), with \( x = \cos \theta \).
		
	\end{definition}
	Note that \( \{ \cos \eta_k \}_{k=1}^n \) denotes a general sequence of nodes. We also note that the operator $\mathcal{G}_n$ also defines an interpolation operator. We refer the sequence of operators $\{\mathcal{G}_n\}$ as Grünwald operators.\par 
	In $1941$, G. Grünwald proved the following convergence theorem for the Chebyshev nodes.
	\begin{theorem}[Grünwald]\cite{grunwald}\label{grun}
		Let $f\in C[-1,1]$, then 
		\[
		\lim\limits_{n\rightarrow\infty}\mathcal{G}_n(f(\cos))(\theta)=f(\cos(\theta)),
		\]
		uniformly on $[0,\pi]$.
	\end{theorem}
	The main difficulty in establishing a convergence theorem arises from the fact that the Lagrange polynomials become unbounded at certain points of $[-1,1]$. Consequently, the following lemma plays a crucial role in proving the result.

	\begin{lemma} [Grünwald] \cite{grunwald}\label{lem}
		\[
		\frac{1}{2}\sum\limits_{k=1}^n \left| P_k\left(\theta - \frac{\pi}{2n}\right) + P_k\left(\theta + \frac{\pi}{2n}\right) \right| < c_1
		\]
		
		where $c_1>0$ is an absolute constant
	\end{lemma} \
	It is interesting to observe that, even though the Lagrange interpolation operators do not converge on the Chebyshev nodes, the sequence of operators \( \mathcal{G}_n \), defined by averaging them, converges uniformly.\par
	By modifying the last step of the proof of Theorem \ref{grun} and using the continuity of $f\in C[0,\pi]$ instead of $f(\cos)$ in \cite{grunwald}, we also have the convergence,
	\[
	\lim\limits_{n\to\infty}\mathcal{G}_n(f)(\theta)=f(\theta)
	\]
	uniformly on $[0,\pi]$.\par
	
	The aim of this article is to investigate for which sets of nodal points the sequence of operators \( \{\mathcal{G}_n\}\) converges. As Grünwald has already resolved this problem in the case of Chebyshev nodes, our objective is to extend this result to other choices of nodal points. This article is structured as follows. In the following section, we address this problem by proving a perturbed version of Theorem~\ref{grun}, thereby identifying a broader class of nodal points for which the convergence result holds. Specifically, we establish the result for the case where the \( \eta_k \)'s in the Definition \ref{Gn} are uniformly distributed and equidistant. We also present an intermediate result that emerges naturally from the proof of our main lemma. Following this, we prove a Voronovskaja-type estimate for the convergence of $\{\mathcal{G}_n\}$. 
	In the subsequent section, we obtain quantitative estimates for the convergence of $\{\mathcal{G}_n\}$ using modulus of continuity.
	
	\section{Convergence on general nodal points}
	In this section, we address the question of whether the convergence result in Theorem \ref{grun} holds for other choices of grid points besides the Chebyshev nodes. Firstly, we note that the convergence of Lagrange interpolation operators critically depends on the selection of grid points. We focus on identifying sets of nodes of the form $\{\cos \eta_k\}_{k=1}^n$ for which Theorem \ref{grun} remains valid. To this end, we establish sufficient conditions under which the convergence stated in Theorem \ref{grun} continues to hold.
	\par
	First, we prove a lemma.
	\begin{lemma}\label{lem1}
		Suppose the nodal points are chosen as $\{\cos \eta_1 , \cos \eta_2 ,\ldots , \cos \eta_n\}$, where
		\[
		\eta_k=\tilde{\theta}_k^{(n)}=\theta_k^{(n)}-\theta_0
		\]
		for $k=1,2,\ldots, n$ and $\theta_0=\theta_0(n)$ is an arbitrary value chosen from the interval $(-\frac{\pi}{2n},\frac{\pi}{2n})$. Then, for $f\in C[-1,1]$ and $\theta\in [0,\pi]$, the Lagrange interpolation operator is given by
		\[
		\mathcal{L}_n(f)(\theta)=\sum\limits_{k=1}^{n}f(\cos\tilde{\theta}_k^{(n)} )l_k(\theta),
		\]
		where
		\[
		l_k(\theta)=\frac{(-1)^{k+1}\cos(n(\theta+\theta_0)).\sin\tilde{\theta}_k^{(n)}}{n(\cos\theta-\cos(\tilde{\theta}_k^{(n)}))}.
		\]
		Then we have,
		\[
		l_k(\theta)=P_k(\theta+\theta_0)S_k(\theta),
		\]
		where $P_k$ denotes the fundamental polynomials of Lagrange interpolation operator on Chebyshev nodes and 
		\[
		S_k(\theta) = \frac{\sin \tilde{\theta}_k^{(n)} \cdot \sin \left( \frac{\theta + \theta_k^{(n)} + \theta_0}{2} \right)}{\sin \theta_k^{(n)} \cdot \sin \left( \frac{\theta + \theta_k^{(n)} - \theta_0}{2} \right)}.
		\]
	\end{lemma}
	\begin{proof}
		Let
		\[
		\zeta(\theta) = \cos n(\theta + \theta_0),
		\]
		where $\theta_0 = \theta_0(n) \in \left(-\frac{\pi}{2n}, \frac{\pi}{2n}\right)$.
		
		Recall that $\theta_k^{(n)} = \frac{(2k - 1)\pi}{2n}$. The roots of $\zeta(\theta)$ are given by $\tilde{\theta}_k^{(n)} = \theta_k^{(n)} - \theta_0.$\par
		Let $x = \cos \theta$. Then,
		\[
		\frac{d}{dx}(\zeta(\theta)) = \frac{d}{d\theta}(\zeta(\theta)) \cdot \frac{d\theta}{dx} 
		= \zeta'(\theta) \cdot \frac{d}{dx}(\arccos x) 
		= \frac{-n \sin(n(\theta + \theta_0))}{- \sin \theta}.
		\]
		Hence,
		\[
		\frac{d}{dx}\zeta(\tilde{\theta}_k^{(n)}) = \frac{(-1)^{k+1} n}{\sin \tilde{\theta}_k^{(n)}}.
		\]
		
		The Lagrange interpolation operator is given by
		\[
		\mathcal{L}_n(f)(\theta) = \sum_{k=1}^n f(x_k) \, l_k(\theta),
		\]
		where
		\[
		l_k(\theta) = \frac{(-1)^{k+1} \cos(n(\theta + \theta_0)) \cdot \sin \tilde{\theta}_k^{(n)}}{n(\cos \theta - \cos \tilde{\theta}_k^{(n)})}.
		\]
		
		Now we simplify $l_k(\theta)$:
		\begin{align*}
			l_k(\theta)
			&= \frac{(-1)^{k+1} \cos(n(\theta + \theta_0)) \cdot \sin \tilde{\theta}_k^{(n)}}{n(\cos \theta - \cos \tilde{\theta}_k^{(n)})} \\
			&= \frac{(-1)^{k+1} \cos n(\theta + \theta_0) \cdot \sin \theta_k^{(n)}}{n(\cos(\theta + \theta_0) - \cos \theta_k^{(n)})} \times \frac{\sin \tilde{\theta}_k^{(n)} \cdot (\cos(\theta + \theta_0) - \cos \theta_k^{(n)})}{\sin \theta_k^{(n)} \cdot (\cos \theta - \cos \tilde{\theta}_k^{(n)})} \\
			&= P_k(\theta + \theta_0) \cdot S_k(\theta),
		\end{align*}
		where $P_k$ denotes the fundamental polynomial of the Lagrange interpolation operator at Chebyshev nodes, and
		\[
		S_k(\theta) = \frac{\sin \tilde{\theta}_k^{(n)} \cdot (\cos(\theta + \theta_0) - \cos \theta_k^{(n)})}{\sin \theta_k^{(n)} \cdot (\cos \theta - \cos \tilde{\theta}_k^{(n)})}.
		\]
		
		Using the identity $\cos a - \cos b = -2 \sin \left( \frac{a + b}{2} \right) \sin \left( \frac{a - b}{2} \right)$ and simplifying, we obtain:
		\[
		S_k(\theta) = \frac{\sin \tilde{\theta}_k^{(n)} \cdot \sin \left( \frac{\theta + \theta_k^{(n)} + \theta_0}{2} \right)}{\sin \theta_k^{(n)} \cdot \sin \left( \frac{\theta + \theta_k^{(n)} - \theta_0}{2} \right)}.
		\]
	\end{proof}
	
	Now, we recall the following theorem established by P. Erdős and G. Grünwald in $1938$. The result gives a uniform bound for the fundamental polynomials of Lagrange interpolation operator on Chebyshev nodes.
	\begin{theorem}[\textbf{P. Erdős, G.  Grünwald}\cite{erdosgrunwald}]\label{erdos}
		In the Lagrange interpolation formula on Chebyshev nodes, 
		\[
		|P_k(\theta)|<\frac{4}{\pi}, 
		\]
		for all $\theta\in\mathbb{R}$, $n$ and $k$.
	\end{theorem} 
	Following this, a minimum bound was obtained by J{\'o}zsef Szabados in $2011$ \cite{szabados}.  Using Theorem \ref{erdos}, we prove the following lemma.
	\begin{lemma}\label{lem2}
		Let $l_k$ be defined in Lemma \ref{lem1}. Then,
		\[
		|l_k(\theta)|<\frac{8}{\pi},
		\]
		for all $k=1,2,\ldots, n$.
	\end{lemma}
	
	\begin{proof}
		By Lemma \ref{lem1}, we have
		$
		l_k(\theta) = P_k(\theta+\theta_0)\, S_k(\theta).$ By Theorem \ref{erdos},
		\[
		|P_k(\theta \pm \theta_0)| < \frac{4}{\pi}, \quad \forall\, \theta \in [0,\pi].
		\]
		Hence, it suffices to show that \( |S_k(\theta)| \leq 2 \) for all \( \theta \in [0,\pi] \), where
		\[
		S_k(\theta) = \frac{\sin(\theta_k^{(n)} - \theta_0) \cdot \sin\left(\frac{\theta + \theta_k^{(n)} + \theta_0}{2}\right)}{\sin \theta_k^{(n)} \cdot \sin\left(\frac{\theta + \theta_k^{(n)} - \theta_0}{2}\right)}.
		\]
		
		We consider the behavior of \( S_k(\theta) \) by analyzing the ranges of \( \theta, \theta_0, \theta_k^{(n)} \) through 10 cases.
		
		\medskip
		
		\noindent
		\textbf{Cases 1--2:} Suppose \( \frac{\theta \pm \theta_0 + \theta_k^{(n)}}{2} \in [0,\frac{\pi}{2}] \) and \( \theta_0 > 0 \).
		
		\begin{itemize}
			\item[(i)] If \( \theta_k^{(n)} \leq \frac{\pi}{2} \), then
			\[
			\sup |S_k(\theta)| \leq S_k(0) = \frac{2\cos \frac{\theta_k^{(n)} - \theta_0}{2} \cdot \sin \frac{\theta_k^{(n)} + \theta_0}{2}}{\sin \theta_k^{(n)}} \leq 2.
			\]
			
			\item[(ii)]If \( \theta_k^{(n)} > \frac{\pi}{2} \), use the identities:
			\[
			\cos \frac{\theta_k^{(n)} - \theta_0}{2} = \sin \frac{\pi - \theta_k^{(n)} + \theta_0}{2}, \quad \sin \theta_k^{(n)} = \sin(\pi - \theta_k^{(n)}),
			\]
			which imply \( |S_k(\theta)| \leq 2 \) by symmetry.
		\end{itemize}
		
		\medskip
		
		\noindent
		\textbf{Cases 3--4:} When \( \theta_0 < 0 \), these follow from Cases 1--2 by replacing \( \theta_0\) by \(-\theta_0 \) and evaluating at \( \theta = \frac{\pi}{2} \), yielding the same bound:
		\[
		|S_k(\theta)| \leq 2.
		\]
		
		\medskip
		
		\noindent
		\textbf{Cases 5--6:} Suppose \( \frac{\theta \pm \theta_0 + \theta_k^{(n)}}{2} \in \left(\frac{\pi}{2}, \pi\right] \) and \( \theta_0 > 0 \). Evaluating at \( \theta = \frac{\pi}{2} \), we get
		\[
		|S_k(\theta)| \leq \frac{\sin(\theta_k^{(n)} - \theta_0)}{\sin \theta_k^{(n)}} \leq \cos \theta_0 + \frac{\cos \theta_k^{(n)} \cdot \sin \theta_0}{\sin \theta_k^{(n)}} \leq 2.
		\]
		
		\medskip
		
		\noindent
		\textbf{Cases 7--8:} When \( \theta_0 < 0 \), proceed as in Cases 5--6 by replacing \( \theta_0\) by \(-\theta_0 \). Again, using standard trigonometric bounds yields \( |S_k(\theta)| \leq 2 \).
		
		\medskip
		
		\noindent
		\textbf{Cases 9--10:} Mixed range cases where one argument lies in \( [0,\frac{\pi}{2}] \) and the other in \( \left(\frac{\pi}{2},\pi\right] \). These transitional cases follow similarly. For instance, evaluating at \( \theta \in \{0, \frac{\pi}{2}, \pi\} \) suffices to show that
		\[
		|S_k(\theta)| \leq 2.
		\]
		
		\medskip
		
		\noindent
		Thus, in all cases, we have \( |S_k(\theta)| \leq 2 \), completing the proof.
	\end{proof}
	
	\begin{table}[H]
		\centering
		\caption{Case-wise bounds on $|S_k(\theta)|$}
		\renewcommand{\arraystretch}{1.5}
		\setlength{\tabcolsep}{10pt}
		\begin{tabular}{|c|c|c|c|c|}
			\hline
			\textbf{Case} & \textbf{$\frac{\theta \pm \theta_0 + \theta_k^{(n)}}{2}$} & \textbf{$\theta_k^{(n)}$} & \textbf{$\theta_0$} & \textbf{Bound on $|S_k(\theta)|$} \\
			\hline
			1 & $[0, \frac{\pi}{2}]$ & $\leq \frac{\pi}{2}$ & $> 0$ & $\leq 2$ \\
			\hline
			2 & $[0, \frac{\pi}{2}]$ & $> \frac{\pi}{2}$ & $> 0$ & $\leq 2$ \\
			\hline
			3 & $[0, \frac{\pi}{2}]$ & $\leq \frac{\pi}{2}$ & $< 0$ & $\leq 2$ \\
			\hline
			4 & $[0, \frac{\pi}{2}]$ & $> \frac{\pi}{2}$ & $< 0$ & $\leq 2$ \\
			\hline
			5 & $(\frac{\pi}{2}, \pi]$ & $\leq \frac{\pi}{2}$ & $> 0$ & $\leq 1$ \\
			\hline
			6 & $(\frac{\pi}{2}, \pi]$ & $> \frac{\pi}{2}$ & $> 0$ & $\leq 2$ \\
			\hline
			7 & $(\frac{\pi}{2}, \pi]$ & $\leq \frac{\pi}{2}$ & $< 0$ & $\leq 2$ \\
			\hline
			8 & $(\frac{\pi}{2}, \pi]$ & $> \frac{\pi}{2}$ & $< 0$ & $\leq 2$ \\
			\hline
			9 & $[0, \frac{\pi}{2}], (\frac{\pi}{2}, \pi]$ & $\leq \frac{\pi}{2}$ & $> 0$ & $\leq 2$ \\
			\hline
			10 & $[0, \frac{\pi}{2}], (\frac{\pi}{2}, \pi]$ & $> \frac{\pi}{2}$ & $> 0$ & $\leq 2$ \\
			\hline
		\end{tabular}
	\end{table}
	
	The following lemma is a perturbed version of Theorem \ref{grun}, providing a broader class of nodal points for which the convergence result holds. We prove this lemma in four steps.
	\begin{lemma}\label{grid}
		Suppose the nodal points are chosen as $\{\cos \eta_1 , \cos \eta_2 ,\ldots , \cos \eta_n\}$, where
		\[
		\eta_k=\tilde{\theta}_k^{(n)}=\theta_k^{(n)}-\theta_0
		\]
		for $k=1,2,\ldots, n$ and $\theta_0=\theta_0(n)$ is an arbitrary value chosen from the interval $(-\frac{\pi}{2n},\frac{\pi}{2n})$. Then, we have
		\[
		\lim\limits_{n\to\infty}\mathcal{G}_n(f(\cos))(\theta)=f(\cos\theta),\ uniformly\ on\ [0,\pi].
		\]
	\end{lemma}
	\begin{proof}
		Let $f\in C[-1,1]$. By Lemma \ref{lem1}, the Lagrange interpolation operator is given by
		\[
		\mathcal{L}_n(f)(\theta)=\sum\limits_{k=1}^nf(\cos\tilde{\theta}_k^{(n)})l_k(\theta),
		\]
		where
		\[
		l_k(\theta)=\frac{(-1)^{k+1}\cos(n(\theta+\theta_0)).\sin\tilde{\theta}_k^{(n)}}{n(\cos\theta-\cos(\tilde{\theta}_k^{(n)})}.
		\]
		Thus, in this case,
		\[
		\mathcal{G}_n(f(\cos))(\theta)=\frac{1}{2}\sum\limits_{k=1}^nf(\cos \tilde\theta_k^{(n)})(l_k(\theta-\frac{\pi}{2n})+l_k(\theta+\frac{\pi}{2n})).
		\]
		To prove this lemma, it suffices to verify the following.
		\begin{equation}\label{bdd}
			\frac{1}{2}\sum\limits_{k=1}^n\left|l_k(\theta-\frac{\pi}{2n})+l_k(\theta+\frac{\pi}{2n})\right |< c_2,
		\end{equation}
		for some $c_2>0$. To prove this, we use Lemma \ref{lem2}.\par
		Let $\theta\neq \tilde{\theta}_k^{(n)}\pm\frac{\pi}{2n}$ or $\theta\neq \frac{k\pi}{n}-\theta_0$, $\frac{(k-1)\pi}{n}-\theta_0$ where $k=1,2,\dots,n$. Following the computational techniques by Grünwald in \cite{grunwald} we have,
		\begin{align*}
			&\frac{1}{2}[l_k(\theta-\frac{\pi}{2n})+l_k(\theta+\frac{\pi}{2n})]\\
			&= \frac{1}{2} \frac{(-1)^{k}\sin\tilde{\theta}_k^{(n)}\sin n(\theta+\theta_0) \sin\theta \sin\frac{\pi}{2n}}{2n} \\
			&\quad \times \frac{1}{\sin\frac{1}{2}(\theta-\theta_0+\theta_k^{(n)}-\frac{\pi}{2n})\sin\frac{1}{2}(\theta+\theta_0-\theta_k^{(n)}-\frac{\pi}{2n})} \\
			&\quad \times \frac{1}{\sin\frac{1}{2}(\theta-\theta_0+\theta_k^{(n)}+\frac{\pi}{2n})\sin\frac{1}{2}(\theta+\theta_0-\theta_k^{(n)}+\frac{\pi}{2n})}.
		\end{align*}
		The equality holds by using the trigonometric identities 
		$\cos (\theta\pm\frac{\pi}{2})=\mp \sin (\theta)$ and $\cos\theta-\cos\eta=-2\sin \frac{1}{2}(\theta+\eta)\sin \frac{1}{2}(\theta-\eta)$ and cross multiplication. We give an estimate for the above term. To obtain an estimate we consider the products $\sin\frac{1}{2}(\theta-\theta_0+\theta_k^{(n)}-\frac{\pi}{2n})\sin\frac{1}{2}(\theta-\theta_0+\theta_k^{(n)}+\frac{\pi}{2n})$ and $\sin\frac{1}{2}(\theta+\theta_0-\theta_k^{(n)}-\frac{\pi}{2n})\sin\frac{1}{2}(\theta+\theta_0-\theta_k^{(n)}+\frac{\pi}{2n})$ and obtain estimates separately.\newline
		We follow the following four steps.\newline
		\textbf{Step 1:}
		Consider,
		\begin{align*}
			&\sin\frac{1}{2}(\theta-\theta_0+\theta_k^{(n)}+\frac{\pi}{2n}) \sin\frac{1}{2}(\theta-\theta_0+\theta_k^{(n)}-\frac{\pi}{2n}) \\
			&= \sin\frac{1}{2}(\theta-\theta_0+\frac{k\pi}{n}) \sin\frac{1}{2}(\theta-\theta_0+\frac{(k-1)\pi}{n}) \\
			&=\sin^2\frac{\theta}{2}\cos(\frac{k\pi}{2n}-\frac{\theta_0}{2})\cos(\frac{(k-1)\pi}{2n}-\frac{\theta_0}{2}) \\
			&+\sin\frac{\theta}{2}\cos\frac{\theta}{2}\cos(\frac{k\pi}{2n}-\frac{\theta_0}{2})\sin(\frac{(k-1)\pi}{2n}-\frac{\theta_0}{2}) \\
			&+\sin\frac{\theta}{2}\cos\frac{\theta}{2}\sin(\frac{k\pi}{2n}-\frac{\theta_0}{2})\cos(\frac{(k-1)\pi}{2n}-\frac{\theta_0}{2}) \\
			&+\cos^2\frac{\theta}{2}\sin(\frac{k\pi}{2n}-\frac{\theta_0}{2})\sin(\frac{(k-1)\pi}{2n}-\frac{\theta_0}{2}).
		\end{align*}
		\textbf{Step 2:} Let $\theta\in[0,\pi]$. Since $\theta_0\in (-\frac{\pi}{2n},\frac{\pi}{2n})$, we have,
		\[
		\frac{1}{2}(\theta+\theta_0-\theta_k^{(n)}+\frac{\pi}{2n})=\frac{1}{2}(\theta+\theta_0-\frac{(k-1)\pi}{n})\in [-\frac{\pi}{2}+\frac{\pi}{4n},\frac{\pi}{2}+\frac{\pi}{4n}]\subset [-\frac{3\pi}{4},\frac{3\pi}{4}]
		\]
		and 
		\[
		\frac{1}{2}(\theta+\theta_0-\theta_k^{(n)}-\frac{\pi}{2n})=\frac{1}{2}(\theta+\theta_0-\frac{k\pi}{n})\in [-\frac{\pi}{2}-\frac{\pi}{4n},\frac{\pi}{2}-\frac{\pi}{4n}]\subset [-\frac{3\pi}{4},\frac{3\pi}{4}].
		\]
		We know that $|\sin x|\geq \frac{4}{3\pi}|x|$ whenever $x\in [-\frac{3\pi}{4},\frac{3\pi}{4}]$. From the above observations we have,
		\[
		|\sin\frac{1}{2}(\theta+\theta_0-\theta_k^{(n)}\underline{+}\frac{\pi}{2n})|\geq \frac{4}{3\pi}\times\frac{1}{2}|\theta+\theta_0-\theta_k^{(n)}\underline{+}\frac{\pi}{2n}|.
		\]
		\textbf{Step 3:}\\
		\textbf{Case 1:} $k\notin \{1 ,n\}$ and $\theta_0\in(-\frac{\pi}{2n},\frac{\pi}{2n})$,
		\[
		\sin^2\frac{\theta}{2}\cos(\frac{k\pi}{2n}-\frac{\theta_0}{2})\cos(\frac{(k-1)\pi}{2n}-\frac{\theta_0}{2})	+ \cos^2\frac{\theta}{2}\sin(\frac{k\pi}{2n}-\frac{\theta_0}{2})\sin(\frac{(k-1)\pi}{2n}-\frac{\theta_0}{2})\geq 0.
		\]
		Thus, in this case we have,
		\[
		\sin \frac{1}{2}(\theta-\theta_0+\theta_k^{(n)}+\frac{\pi}{2n}) \sin\frac{1}{2}(\theta-\theta_0+\theta_k^{(n)}-\frac{\pi}{2n})\geq \frac{1}{2}\sin\theta\sin\tilde\theta_k^{(n)}
		\]
		and hence,
		\begin{equation}
			\frac{1}{2} \left| \frac{ \sin\theta \, \sin\tilde{\theta}_k^{(n)} }{ \sin\left( \frac{1}{2}(\theta + \tilde\theta_k^{(n)} + \frac{\pi}{2n}) \right) \sin\left( \frac{1}{2}(\theta + \tilde\theta_k^{(n)} - \frac{\pi}{2n}) \right) } \right| \leq 1.
		\end{equation}
		Following the proof of Grünwald, for $k\neq 1,n$, $\theta\neq \tilde\theta_k^{(n)}\underline{+}\frac{\pi}{2n}$ and $\theta_0\in(-\frac{\pi}{2n},\frac{\pi}{2n})$, we have,\par
		
		$\frac{1}{2} \left| l_k\left(\theta + \frac{\pi}{2n}\right) + l_k\left(\theta - \frac{\pi}{2n}\right) \right|$
		\begin{align}\label{ine1}
			&\leq \left| \frac{\sin\frac{\pi}{2n}}{2n} \right| 
			\left| \frac{1}{
				\sin\left( \frac{1}{2}\left( \theta - \tilde\theta_k^{(n)} + \frac{\pi}{2n} \right) \right)
				\sin\left( \frac{1}{2}\left( \theta - \tilde\theta_k^{(n)} - \frac{\pi}{2n} \right) \right)
			} \right|  \nonumber \\
			&\leq \frac{9\pi^3}{64n^2} \cdot 
			\frac{1}{
				\left| \theta - \tilde\theta_k^{(n)} + \frac{\pi}{2n} \right|
				\left| \theta - \tilde\theta_k^{(n)} - \frac{\pi}{2n} \right|
			} \nonumber  \\
			&\leq \frac{9\pi^3}{64n^2} \cdot 
			\frac{1}{\left( \theta - \tilde\theta_k^{(n)} - \frac{\pi}{2n} \right)^2}
		\end{align}
		\textbf{Case 2:} $k\in\{1,n\}$ and $\theta_0\in(-\frac{\pi}{2n},\frac{\pi}{2n})$.
		Let $\delta>0$ be fixed.
		Suppose $|\theta-\tilde\theta_1^{(n)}|>\delta$ or $|\theta-\tilde\theta_n^{(n)}|>\delta$, for sufficiently large $n$, then we have,
		\begin{equation}
			\frac{1}{2} \left| \frac{\sin \theta \, \sin \tilde{\theta}_k^{(n)}}{
				\sin \frac{1}{2} \left( \theta + \tilde{\theta}_k^{(n)} + \frac{\pi}{2n} \right) 
				\cdot 
				\sin \frac{1}{2} \left( \theta + \tilde{\theta}_k^{(n)} - \frac{\pi}{2n} \right) 
			} \right| \leq C,
		\end{equation}
		where $C$ is a constant. Therefore, in this case,
		\begin{equation}\label{ine2}
			\frac{1}{2} \left| l_k\left(\theta + \frac{\pi}{2n}\right) + l_k\left(\theta - \frac{\pi}{2n}\right) \right|
			\leq C \frac{9 \pi^3}{64 n^2} \frac{1}{\left(\theta - \tilde{\theta}_k^{(n)} - \frac{\pi}{2n}\right)^2}.
		\end{equation}
		In both cases above, the inequalities \ref{ine1} and \ref{ine2} holds trivially for $\theta=0,\pi$.\\
		\textbf{Step 4:}\\
		Let $\theta$ be fixed and 
		\[
		1>\tilde x_1>\tilde x_2>...>\tilde x_j>x=\cos\theta>\tilde x_{j+1}>...>\tilde x_n>-1,
		\]
		where $\tilde x_k=\cos \tilde \theta_k^{(n)}$.
		By Lemma \ref{lem2}, we have $|l_k(\theta)|\leq \frac{8}{\pi}$ for all $\theta\in [0,\pi]$, $k$ and $n$. Now,\newline
		$\frac{1}{2}\sum\limits_{k=1}^{n}|l_k(\theta-\frac{\pi}{2n})+l_k(\theta+\frac{\pi}{2n})|$
		\begin{align}\label{equat2}
			&\leq \frac{48}{\pi} + \frac{1}{2} \sum_{\substack{2 \leq k \leq n - 1 \\ k \neq j-2,\, j-1,\, j,\, j+1}} 
			\left| l_k\left( \theta - \frac{\pi}{2n} \right) + l_k\left( \theta + \frac{\pi}{2n} \right) \right| \nonumber \\
			&\leq \frac{48}{\pi} + \sum_{\substack{2 \leq k \leq n - 1 \\ k \neq j-2,\, j-1,\, j,\, j+1}} 
			\frac{9\pi^3}{64n^2} \cdot \frac{1}{\left( \theta - \tilde{\theta}_k^{(n)} - \frac{\pi}{2n} \right)^2} \nonumber  \\
			&\leq \frac{48}{\pi} + \frac{9\pi^3}{64n^2} \sum_{k=1}^{\infty} 
			\left( \frac{2n}{\pi} \right)^2 \cdot \frac{1}{k^2} = c_2.
		\end{align}
		We use,
		\[
		|\theta - \tilde{\theta}_k^{(n)} - \frac{\pi}{2n}| \geq |\tilde{\theta}_k^{(n)} - \tilde{\theta}_{j+1}^{(n)}| = \left| \frac{(k - j - 1)\pi}{2n} \right|
		\]
		or
		\[
		|\theta - \tilde{\theta}_k^{(n)} - \frac{\pi}{2n}| \geq |\tilde{\theta}_k^{(n)} - \tilde{\theta}_{j-1}^{(n)}| = \left| \frac{(k - j + 1)\pi}{2n} \right|.
		\]
		Hence, inequality~\eqref{bdd} is verified.
		\par
		
		From \ref{ine2} and \ref{equat2},  for sufficiently large $n$ and for a fixed $\delta>0$, we have,
		\begin{equation}\label{equat3}
			\sum_{\substack{1 \leq k \leq n \\ |\theta - \tilde{\theta}_k^{(n)}| > \delta}}
			\frac{1}{2} \left| l_k\left(\theta - \frac{\pi}{2n}\right) + l_k\left(\theta + \frac{\pi}{2n}\right) \right| 
			= \mathcal{O}\left(\frac{1}{n}\right).
		\end{equation}

		Let $\epsilon>0$, applying the continuity of the function $f(\cos(.))$, there exists a $\delta>0$ such that 
		\[
		|f(\cos\theta)-f(\cos\tilde\theta_k^{(n)})|< \epsilon
		\]
		whenever $|\theta-\tilde\theta_k^{(n)}|\leq \delta$.
		Now, we have,
		\begin{multline}\label{equat4}
			\left| \frac{1}{2}\left\{ \mathcal{L}_n(f)\left(\theta - \frac{\pi}{2n}\right) + \mathcal{L}_n(f)\left(\theta + \frac{\pi}{2n}\right) \right\} - f(\cos\theta) \right| \leq \\
			\frac{1}{2} \sum_{\substack{1 \leq k \leq n \\ |\theta - \tilde{\theta}_k^{(n)}| \leq \delta}} 
			| f(\cos\theta) - f(\cos \tilde{\theta}_k^{(n)}) | \, 
			\left| l_k\left(\theta - \frac{\pi}{2n}\right) + l_k\left(\theta + \frac{\pi}{2n}\right) \right| + \\
			\frac{1}{2} \sum_{\substack{1 \leq k \leq n \\ |\theta - \tilde{\theta}_k^{(n)}| > \delta}} 
			| f(\cos\theta) - f(\cos \tilde{\theta}_k^{(n)}) | \, 
			\left| l_k\left(\theta - \frac{\pi}{2n}\right) + l_k\left(\theta + \frac{\pi}{2n}\right) \right|.
		\end{multline}
		
		The first expression in the above inequality can be made sufficiently small using the continuity of the function $f(\cos)$ and the inequality \ref{equat2}. The second expression can be made sufficiently small for large $n$ using the boundedness of the function $f(\cos)$ on $[0,\pi]$ and \ref{equat3}.
	\end{proof}
	In the Lemma \ref{grid}, if we choose $\theta_0=0$, then, Theorem \ref{grun} follows. In this case, $l_k(\theta)=P_k(\theta)$, for $\theta\in [0,\pi]$ and $k=1,2,\ldots n$. In what follows, our aim is to characterize the set of points $\{\eta_0,\eta_1,\ldots,\eta_{n-1}\}$ in $[0,\pi]$ for which Theorem \ref{grun} holds. Although a complete characterization of such points is not achieved, we establish a result that generalizes Lemma \ref{grun}. Specifically, we establish a convergence result for the case where the given points are equidistant and uniformly distributed.\par
	\begin{theorem}\label{grid1}
		Let $n\geq 2$ and $\eta_0<\eta_1<\ldots <\eta_{n-1}$ be equidistant points in $[0,\pi]$ such that 
		\begin{enumerate}
			\item $\eta_0\leq \frac{\pi}{n}$.
			\item $\eta_{i+1}-\eta_{i}=\frac{\pi}{n}-\beta_n$ for some $\beta_n$ such that $0\leq \beta_n<\frac{\pi}{n}$ and $i=0,1,\ldots,n-2$.
		\end{enumerate}
		Then, for $f\in C[-1,1]$, the following convergence holds with respect to the nodes:
		\[
		\{\cos\eta_0,\cos(\eta_1+\beta_n),\ldots,\cos(\eta_{n-1}+(n-1)\beta_n)\},
		\]
		\[
		\lim\limits_{n\to\infty}\mathcal{G}_n(f(\cos))(\theta)=f(\cos\theta)
		\]
		uniformly on $[0,\pi]$. 
	\end{theorem}
	
	\begin{proof}
		We have, $\eta_{i+1}-\eta_{i}=\frac{\pi}{n}-\beta_n$, where $0\leq \beta_n\leq \frac{\pi}{n}$.
		Thus, 
		\begin{flalign*}
			&\eta_0 = (\eta_0-\frac{\pi}{2n})+\theta_1^{(n)}, \\
			&\eta_1+\beta_n = \eta_0+\frac{\pi}{n} = (\eta_0-\frac{\pi}{2n})+\theta_2^{(n)}, \\
			&\eta_2+2\beta_n = \eta_0+\frac{2\pi}{n} = (\eta_0-\frac{\pi}{2n})+\theta_3^{(n)}, \\
			&\vdots \\
			&\eta_{n-1}+(n-1)\beta_n = \eta_{0}+\frac{(n-1)\pi}{n} = (\eta_0-\frac{\pi}{2n})+\theta_n^{(n)}.
		\end{flalign*}
		We have, $\eta_0\leq \frac{\pi}{n}$, so that $\eta_0-\frac{\pi}{2n}\in(-\frac{\pi}{2n},\frac{\pi}{2n})$. 
		Now, it is easy to apply Lemma \ref{grid} by choosing $\theta_0=\eta_0-\frac{\pi}{2n}$. Hence the result follows.
	\end{proof}
	
	\begin{corollary}
		Let $n\geq 2$ and $\eta_0<\eta_1<\ldots <\eta_{n-1}$ be equidistant points in $[0,\pi]$ such that 
		\begin{enumerate}
			\item $\eta_0\leq \frac{\pi}{n}$.
			\item $\eta_{i+1}-\eta_{i}=\frac{s}{n}\ for\ some\ s\ such\ that\ 0<s\leq \pi \ and\ i=0,1,\ldots,n-2$.
		\end{enumerate}
		Then, for $f\in C[-1,1]$, the following convergence holds with respect to the nodes:
		\[ \{\cos\eta_0,\cos(\eta_1+\frac{\pi-s}{n}),\ldots,\cos(\eta_{n-1}+\frac{(n-1)(\pi-s)}{n})\}:
		\]
		\[
		\lim\limits_{n\to\infty}\mathcal{G}_n(f(\cos))(\theta)=f(\cos\theta)
		\]
		uniformly on $[0,\pi]$. 
	\end{corollary}
	\begin{proof}
		We have, $\eta_{i+1}-\eta_{i}=\frac{s}{n}=\frac{\pi}{n}-\frac{\pi-s}{n}$ for $i=0,1,\ldots, n$. Let $\beta_n=\frac{\pi-s}{n}$ and the result follows by applying Theorem \ref{grid1}.
	\end{proof}
	We recall the definition of uniform distribution of grid points on an interval $I$.
	
	\begin{definition}\label{uniform}
		A sequence of grids $\{S_n=\{x_i^{(n)}:i=0,1,...,n\}\}_n$ belonging to an interval $I$ is uniformly distributed if 
		\[
		\sum\limits_{i=1}^n|\frac{|I|}{n}-(x_i^{(n)}-x_{i-1}^{(n)})|=\mathcal{O}(\frac{1}{n})\ (n\rightarrow\infty)
		\]
		with $|I|$ being the length of $I$.
	\end{definition}
	\begin{remark}
		From Lemma \ref{grid}, it is evident that convergence holds for the nodal points, $\{\cos\tilde\theta_1^{(n)}, \cos\tilde\theta_2^{(n)}, \ldots, \cos\tilde\theta_n^{(n)}\}$ where $\tilde\theta_k^{(n)}=\theta_k^{(n)}+\frac{\pi}{n^l}$, $l\geq2,\ l\in\mathbb{N}$ which is uniformly distributed.
	\end{remark}\par 
	We have the following result.
	\begin{theorem}
		Let $n\geq 2$ and $\eta_0<\eta_1<\ldots <\eta_{n-1}$ be uniformly distributed equidistant points in $[0,\pi]$ such that $\eta_0\leq \frac{\pi}{n}$. Then there exist $\beta_n$ satisfying $0\leq \beta_n< \frac{\pi}{n}$ such that, for $f\in C[-1,1]$, the following convergence holds with respect to the nodes:
		\[ \{\cos\eta_0,\cos(\eta_1+\beta_n),\ldots,\cos(\eta_{n-1}+(n-1)\beta_n)\},
		\]
		\[
		\lim\limits_{n\to\infty}\mathcal{G}_n(f(\cos))(\theta)=f(\cos\theta),
		\]
		uniformly on $[0,\pi]$. 
	\end{theorem}
	\begin{proof}
		Since  $\eta_0<\eta_1<\ldots <\eta_{n-1}$ equidistant, we have, $\eta_{i+1}-\eta_i=\alpha_n$ for $i=1,2,\ldots,n-2$. Since it is also uniformly distributed we have 
		\[
		\sum\limits_{i=1}^{n-1}|\frac{\pi}{n}-(\eta_i-\eta_{i-1})|=\mathcal{O}(\frac{1}{n-1})\ (n\rightarrow\infty).
		\]
		Now, $0<\alpha_n\leq \frac{\pi}{n}$, thus, 
		\[
		\frac{\pi}{n}-\alpha_n=\mathcal{O}(\frac{1}{(n-1)^2})(n\to\infty)
		\]
		Set $\beta_n=\frac{\pi}{n}-\alpha_n$. Therefore, $\alpha_n=\frac{\pi}{n}-\beta_n$. Clearly, $\beta_n< \frac{\pi}{n}$. Applying Theorem \ref{grid1}, we get the result.
	\end{proof}
	The following is an intermediate result that can be derived from the proof of Lemma \ref{grid}.
	\begin{corollary}
		For $n=1,2,\ldots$, let $s_1^{(n)} <s_2^{(n)}< \ldots < s_n^{(n)}$ be points such that $s_i^{(n)}\in (-\frac{\pi}{2n},\frac{\pi}{2n})$ and $\eta_k^{(n)}=\theta_k^{(n)}+s_k^{(n)}$ for $k=1,2,\ldots n$ and $i=1,2,\ldots$. Define 
		\[
		\tilde{\mathcal{G}}_n(f)(\theta):= \frac{1}{2}\sum \limits_{k=1}^n f(\cos \eta_k^{(n)})\{\tilde{P}_k(\theta+\frac{\pi}{2n})+\tilde{P}_k(\theta-\frac{\pi}{2n})\},
		\]
		where 
		\[
		\tilde{P}_k(\theta)=\frac{(-1)^{k+1}\cos(n(\theta+s_k^{(n)})).\sin {\eta}_k^{(n)}}{n(\cos\theta-\cos({\eta}_k^{(n)})}.
		\]
		Then we have 
		\[
		\lim\limits_{n\to\infty}\tilde{\mathcal{G}}_n(f(\cos))(\theta)=f(\cos\theta).
		\]
		in $C[0,\pi]$.
	\end{corollary}
	\begin{proof}
		This follows from the proof of Lemma \ref{grid}.
	\end{proof}
	
	In the following section, we obtain a Voronovskaja-type estimate for the convergence of the sequence $\{\mathcal G_n\}$. Thereafter, we derive quantitative estimates for the rate of convergence of $\{\mathcal G_n\}$ using the first modulus of smoothness. Throughout this analysis, the nodes are taken to be the perturbed nodes $\{\cos \tilde{\theta}_k^{(n)}\}_{k=1}^n$ as described in Lemma~\ref{grid}.
	
	\section{A Voronovskaja-type estimate}
	Voronovskaja-type theorems form an important class of results in approximation theory, describing the asymptotic behaviour of sequences of positive linear operators that approximate continuous functions. The classical result was first established by E.~Voronovskaja in 1932 for the Bernstein polynomials $\{B_n\}$ \cite{voronovskaja}. For every $f \in C^2[0,1]$ ($f$ is twice differentiable and the second derivative is continuous), we have
	\[
	\lim_{n \to \infty} n \big(B_n(f;x) - f(x)\big) = \frac{1}{2} x(1-x) f''(x). \qquad x \in [0,1].
	\]
	This theorem provided not only the convergence of $B_n(f;x)$ to $f(x)$ but also the rate and the precise form of this convergence. It established a fundamental connection between the approximation process and the smoothness of the function through its derivatives. Since then, Voronovskaja’s result has been extended to various classes of operators. For recent results, we refer to \cite{acar, agrawal, popa}. For related results concerning sequences of non-positive operators, see \cite{bianca}.
	
	\par
	
	We now prove a Voronovskaja-type estimate for the sequence $\{\mathcal{G}_n\}$. A full Voronovskaja-type theorem cannot be obtained in this setting; nevertheless, we derive an estimate for the convergence.
	For this purpose, we obtain the following.
\begin{proposition}\label{ineq}
	Let $\{\mathcal{G}_n\}$ be defined on the perturbed nodes $\{\cos\theta_k^{(n)}\}$. Let $|\tilde\theta_k^{n}-\theta|>\kappa_n$ for $k=1,2,\ldots, n$ for some $\kappa_n>\frac{\pi}{2n}$. Then, 
	\begin{equation}\label{ineq1}
		\frac{1}{2}\sum\limits_{k=1}^{n}|l_k(\theta-\frac{\pi}{2n})+l_k(\theta+\frac{\pi}{2n})|\leq\frac{9\pi^3}{64n^2} \cdot 
		\sum\limits_{k=1}^{n}\frac{1}{\left( \theta - \tilde\theta_k^{(n)} - \frac{\pi}{2n} \right)^2}.
	\end{equation}
\end{proposition}

\begin{proof}
	It suffices to prove the inequality \ref{ineq1} for $k=1,n$, since this holds for all $k=2,\ldots, n-1$ by \ref{ine1}. To prove this, it is enough to verify the following inequality for $k=1,n$.
	\begin{equation}\label{ine3}
		\frac{1}{2} \left| \frac{\sin \theta \, \sin \tilde{\theta}_k^{(n)}}{
			\sin \frac{1}{2} \left( \theta + \tilde{\theta}_k^{(n)} + \frac{\pi}{2n} \right) 
			\cdot 
			\sin \frac{1}{2} \left( \theta + \tilde{\theta}_k^{(n)} - \frac{\pi}{2n} \right) } \right| \leq 1,
	\end{equation}
	
	\textbf{Case 1: k=1} $|\theta-\tilde\theta_1^{(n)}|>\kappa_n$. Then either 
	$\theta>\frac{\pi}{2n}-\theta_0+\kappa_n$ or 
	$\theta<\frac{\pi}{2n}-\theta_0-\kappa_n$. Since $\kappa_n>\frac{\pi}{2n}$ and 
	$\theta\in [0,\pi]$, the second condition cannot be true. Therefore
	\[
	\theta>\frac{\pi}{2n}-\theta_0+\kappa_n.
	\]
	
	We have 
	\begin{equation*}
		\begin{split}
			\frac{\sin \theta \, \sin \tilde{\theta}_1^{(n)}}{
				\sin \frac{1}{2} \left( \theta + \tilde{\theta}_1^{(n)} + \frac{\pi}{2n} \right) 
				\cdot 
				\sin \frac{1}{2} \left( \theta + \tilde{\theta}_1^{(n)} - \frac{\pi}{2n} \right) }
			&= \frac{\sin \theta \, \sin \left(\frac{\pi}{2n}-\theta_0\right)}{
				\sin \frac{1}{2} \left( \theta + \frac{\pi}{n}-\theta_0 \right) 
				\cdot 
				\sin \frac{1}{2}(\theta-\theta_0)}\\[2mm]
			&= \frac{2\cos \frac{\theta}{2}\sin \frac{\theta}{2} \, 
				\sin \left(\frac{\pi}{2n}-\theta_0\right)}{
				\sin \frac{1}{2} \left( \theta + \frac{\pi}{n}-\theta_0 \right)
				\sin \frac{1}{2}(\theta-\theta_0)}.
		\end{split}
	\end{equation*}
	
	For $\theta\leq \frac{\pi}{2}$, we have 
	\[
	\sin \frac{\theta}{2}\leq 
	\sin \frac{1}{2} \left( \theta + \frac{\pi}{n}-\theta_0 \right), \quad
	\text{and} \quad
	\sin \left(\frac{\pi}{2n}-\theta_0\right)
	\leq 
	\sin \frac{1}{2}(\theta-\theta_0),
	\]
	since $\kappa_n>\frac{\pi}{2n}$.
	
	For $\theta>\frac{\pi}{2}$, we have
	\[
	\cos\frac{\theta}{2}\leq 
	\sin \frac{1}{2} \left( \theta + \frac{\pi}{n}-\theta_0 \right)
	= \cos\frac{1}{2} \left( \pi-\theta -\frac{\pi}{n}+\theta_0 \right),
	\]
	since $\cos\theta$ is non-increasing. 
	
	Note that 
	\[
	\frac{\theta}{2}\geq 
	\frac{1}{2} \left(\pi-\theta -\frac{\pi}{n}+\theta_0\right),
	\]
	which holds since
	\[
	\theta> \frac{\pi}{2}\geq
	\frac{\pi}{2}-\frac{\pi}{2n}+\frac{\theta_0}{2}.
	\]
	In either case, inequality \ref{ine3} holds.
	
	\textbf{Case 2: k=n} $|\theta-\tilde\theta_n^{(n)}|>\kappa_n$.
	
	In this case, either 
	\[
	\theta>\frac{(2n-1)\pi}{2n}-\theta_0+\kappa_n \quad \text{or} \quad 
	\theta<\frac{(2n-1)\pi}{2n}-\theta_0-\kappa_n.
	\]
	Since $\kappa_n>\frac{\pi}{2n}$, the first case cannot occur. Hence, 
	\[
	\theta<\frac{(2n-1)\pi}{2n}-\theta_0-\kappa_n = \pi-\frac{\pi}{n}-\theta_0-\kappa_n.
	\]
	
	We have 
	\begin{equation*}
		\begin{split}
			\frac{\sin \theta \, \sin \tilde{\theta}_n^{(n)}}{
				\sin \frac{1}{2} \left( \theta + \tilde{\theta}_n^{(n)} + \frac{\pi}{2n} \right) 
				\cdot 
				\sin \frac{1}{2} \left( \theta + \tilde{\theta}_n^{(n)} - \frac{\pi}{2n} \right) } 
			&= \frac{\sin \theta \, \sin (\frac{\pi}{n}-\theta_0)}{
				\cos \frac{1}{2} \left( \theta-\theta_0 \right) 
				\cdot 
				\cos \frac{1}{2}(\theta-\frac{\pi}{n}-\theta_0)}\\
			&= \frac{2\sin \frac{\theta}{2}\cos \frac{\theta}{2} \,\sin (\frac{\pi}{n}-\theta_0)}{
				\cos \frac{1}{2} \left( \theta-\theta_0 \right) 
				\cdot 
				\cos \frac{1}{2}(\theta-\frac{\pi}{n}-\theta_0)}.
		\end{split}
	\end{equation*}
	
	For all $\theta\in[0,\pi]$, since $\cos\theta$ is non-increasing,
	\[
	\cos\frac{\theta}{2}\leq \cos \frac{1}{2}(\theta-\theta_0).
	\]
	
	Let $\theta\leq \frac{\pi}{2}$. Since $\cos\theta$ is increasing on $[0,\frac{\pi}{2}]$, we have 
	\[
	\sin\frac{\theta}{2} = \cos \Big(\frac{\pi}{2}-\frac{\theta}{2}\Big) \leq \cos \frac{1}{2}\big(\theta-\frac{\pi}{n}-\theta_0\big),
	\]
	since 
	\[
	\frac{\pi}{2}-\frac{\theta}{2} > \frac{1}{2}(\theta-\frac{\pi}{n}-\theta_0).
	\]
	
	For $\theta>\frac{\pi}{2}$, we have 
	\[
	\sin\frac{\theta}{2} = \cos\Big(\frac{\pi}{2}-\frac{\theta}{2}\Big).
	\]
	Since $\cos\theta$ is non-increasing on $[0,\pi]$,
	\[
	\cos\Big(\frac{\pi}{2}-\frac{\theta}{2}\Big) \leq \cos \frac{1}{2}(\theta-\frac{\pi}{n}-\theta_0).
	\]
	
	Now, 
	\[
	\frac{\pi}{2}-\frac{\theta}{2} \geq \frac{1}{2}(\theta-\frac{\pi}{n}-\theta_0),
	\]
	which holds if and only if 
	$
	\theta < 3\theta_0+\pi-\frac{\pi}{n},
	$ 
	which is true since 
	$
	\theta<\pi-\frac{\pi}{n}-\theta_0-\kappa_n.
	$
	
	Hence, by the above cases, \ref{ine3} is verified.
\end{proof}
Let $c_2$ be the constant in the equation \ref{equat2}.
	\begin{theorem}\label{voronov}
		Let $f\in C^2[0,\pi]$ and $\theta\in[0,\pi]$. Then
		\[
		\limsup_{n\to\infty} n^{1/3}\, |\mathcal{G}_n(f)(\theta)-f(\theta)|
		\le 
		\left(\frac{9\pi^4}{64}+c_2\right)|f'(\theta)|
		+\frac{9\pi^5}{128}|f''(\theta)|
		+\frac{9M\pi^5}{64},
		\]
		where $M=\|f''\|_\infty$.
	\end{theorem}
	
	\begin{proof}
		Let $f\in C^2[0,\pi]$ and $t,\theta\in[0,\pi]$.  
		By Taylor's theorem,
		\[
		f(t)-f(\theta)
		=
		f'(\theta)(t-\theta)
		+\frac{f''(\theta)}{2}(t-\theta)^2
		+\frac{f''(c)-f''(\theta)}{2}(t-\theta)^2,
		\]
		where $c$ lies between $t$ and $\theta$.
		
		From equation \ref{equat2}, we have 
		\begin{equation}\label{ine4}
			\frac12\sum_{k=1}^n
			\bigl|
			l_k(\theta-\tfrac{\pi}{2n})
			+
			l_k(\theta+\tfrac{\pi}{2n})
			\bigr|
			\le c_2.
		\end{equation}
		
		Let $e_0(t)=t$ for $t\in[0,\pi]$ and define $\epsilon_n=n^{-1/3}$.
		Then
		\begin{align*}
			|\mathcal{G}_n(e_0)(\theta)-e_0(\theta)|
			&\le
			\frac12
			\sum_{\substack{1\le k\le n\\ |\theta-\tilde\theta_k^{(n)}|<\frac{\pi}{2n}+\epsilon_n}}
			|\tilde\theta_k^{(n)}-\theta|
			\bigl|
			l_k(\theta-\tfrac{\pi}{2n})+l_k(\theta+\tfrac{\pi}{2n})
			\bigr|
			\\
			&\quad+
			\frac12
			\sum_{\substack{1\le k\le n\\ |\theta-\tilde\theta_k^{(n)}|\ge \frac{\pi}{2n}+\epsilon_n}}
			|\tilde\theta_k^{(n)}-\theta|
			\bigl|
			l_k(\theta-\tfrac{\pi}{2n})+l_k(\theta+\tfrac{\pi}{2n})
			\bigr|
			\\[2mm]
			&\le
			\left(\frac{\pi}{2n}+\epsilon_n\right)c_2
			+
			\pi\cdot \frac{9\pi^3}{64n^2}
			\sum_{k=1}^n \frac{1}{(\theta-\tilde\theta_k^{(n)}-\frac{\pi}{2n})^2}
			\\
			&\le
			\left(\frac{\pi}{2n}+\epsilon_n\right)c_2
			+
			\frac{9\pi^4}{64}\frac{1}{n\epsilon_n^2},
		\end{align*}
		using \eqref{ine4}, $|\theta-\tilde\theta_k^{(n)}|\le \pi$, and Proposition~\ref{ineq} with 
		$\kappa_n=\frac{\pi}{2n}+\epsilon_n>\frac{\pi}{2n}$.
		
		Similarly,
		\begin{align*}
			|\mathcal{G}_n((e_0-\theta)^2)(\theta)|
			&\le
			\left(\frac{\pi}{2n}+\epsilon_n\right)^2 c_2
			+
			\frac{9\pi^5}{64}\frac{1}{n\epsilon_n^2}
			\\
			&=
			\left(\frac{\pi^2}{4n^2}+\epsilon_n^2+\frac{\pi\epsilon_n}{n}\right)c_2
			+\frac{9\pi^5}{64}\frac{1}{n\epsilon_n^2}.
		\end{align*}
		
		Let $M=\|f''\|_\infty$. Then
		\[
		|\mathcal{G}_n\big((e_0-\theta)^2(f''(c)-f''(\theta))\big)(\theta)|
		\le
		2M\left\{
		\left(\frac{\pi^2}{4n^2}+\epsilon_n^2+\frac{\pi\epsilon_n}{n}\right)c_2
		+\frac{9\pi^5}{64}\frac{1}{n\epsilon_n^2}
		\right\}.
		\]
		
		Now multiply these inequalities by $n\epsilon_n^2$:
		
		\[
		n\epsilon_n^2|\mathcal{G}_n(e_0)(\theta)-e_0(\theta)|
		\le
		\left(\frac{\pi\epsilon_n^2}{2}+n\epsilon_n^3\right)c_2
		+\frac{9\pi^4}{64},
		\]
		
		\[
		n\epsilon_n^2|\mathcal{G}_n((e_0-\theta)^2)(\theta)|
		\le
		\left(
		\frac{\pi^2\epsilon_n^2}{4n}
		+n\epsilon_n^4
		+\pi\epsilon_n^3
		\right)c_2
		+
		\frac{9\pi^5}{64},
		\]
		
		\[
		n\epsilon_n^2
		\big|
		\mathcal{G}_n((e_0-\theta)^2(f''(c)-f''(\theta)))(\theta)
		\big|
		\le
		2M
		\left(
		\frac{\pi^2\epsilon_n^2}{4n}
		+n\epsilon_n^4
		+\pi\epsilon_n^3
		\right)c_2
		+
		\frac{18M\pi^5}{64}.
		\]
		
		Since $\epsilon_n=n^{-1/3}$, we have
		\[
		\epsilon_n\to 0,\qquad
		n\epsilon_n^4\to 0,\qquad
		\frac{1}{n\epsilon_n^2}\to 0,\quad as\ n\to\infty, \qquad
		n\epsilon_n^3=1.
		\]
		
		Combining the estimates from the Taylor expansion shows that
		\[
		\limsup_{n\to\infty} n^{1/3}\,|\mathcal{G}_n(f)(\theta)-f(\theta)|
		\le
		\left(\frac{9\pi^4}{64}+c_2\right)|f'(\theta)|
		+\frac{9\pi^5}{128}|f''(\theta)|
		+\frac{9M\pi^5}{64},
		\]
		as claimed.
	\end{proof}
	
	\section{Quantitative estimates}
	In this section, we obtain a quantitative estimate for the convergence of the sequence of operators $\{\mathcal{G}_n\}$ using modulus of continuity. Such quantitative convergence results have been extensively studied for positive linear operators; see, for example, \cite{devore, ditzian, swetits}.
	
	We recall the definition of the classical modulus of continuity.
	
\begin{definition}
	Let $f:[a,b]\to \mathbb{R}$ be a bounded function.  
	The first-order modulus of smoothness (or modulus of continuity) of $f$ with argument $\delta>0$ is defined by  
	\[
	\omega_1(f,\delta)
	= \sup\{\, |f(x)-f(y)| : |x-y|\le \delta,\; x,y\in[a,b] \,\}.
	\]
\end{definition}

In the following, we obtain a quantitative estimate for $\{\mathcal{G}_n\}$ using modulus of continuity. Let $\|\cdot\|_\infty$ denotes the supremum norm on $[0,\pi]$.

\begin{theorem}
	Let $f\in C[0,\pi]$. Then 
	\[
	\|\mathcal G_n(f)-f\|_\infty
	\leq 
	(c_2+1)\omega(f,\mu_n),
	\]
where 
\[
\mu_n:=\left(\frac{\pi}{2n}+\frac{1}{n^{\frac{1}{3}}}\right)c_2
+
\frac{9\pi^4}{64}\frac{1}{n^{\frac{1}{3}}}.
\]
\end{theorem}

\begin{proof}
Let $f\in C[0,\pi]$, $t,\theta\in [0,\pi]$ and $\delta>0$, then by the properties of the modulus of continuity,
\begin{equation*}
	\begin{split}
|f(t)-f(\theta)|&\leq \omega(f,|t-\theta|)\\
&\leq (1+|t-\theta|\delta^{-1})\omega(f,\delta)
\end{split}
\end{equation*}
Now, 
\begin{equation*}
	\begin{split}
\mathcal{G}_n(f)(\theta)-f(\theta)&=\frac{1}{2}\sum\limits_{k=1}^n(f(\theta_k^{(n)})-f(\theta))\Big[	l_k(\theta-\tfrac{\pi}{2n})
+
l_k(\theta+\tfrac{\pi}{2n})\Big ].
\end{split}
\end{equation*}
Thus 
\[
\begin{aligned}
	|\mathcal{G}_n(f)(\theta)-f(\theta)|
	&\le \frac{1}{2}\sum_{k=1}^n 
	\bigl| f(\theta_k^{(n)}) - f(\theta)\bigr|
	\bigl|\,l_k(\theta-\tfrac{\pi}{2n}) +
	l_k(\theta+\tfrac{\pi}{2n}) \bigr| \\[6pt]
	&\le \frac{1}{2}\sum_{k=1}^n 
	\Bigl|l_k(\theta-\tfrac{\pi}{2n}) +
	l_k(\theta+\tfrac{\pi}{2n}) \Bigr|
	\,\omega(f,\delta) \\[6pt]
	&\quad + \frac{1}{2}\sum_{k=1}^n 
	\bigl|\theta_k^{(n)} - \theta \bigr|
	\Bigl|l_k(\theta-\tfrac{\pi}{2n}) +
	l_k(\theta+\tfrac{\pi}{2n}) \Bigr|\delta^{-1}
	\,\omega(f,\delta) \\[6pt]
	&\le c_2\,\omega(f,\delta)
	+ \frac{1}{2}\sum_{k=1}^n 
	\bigl|\theta_k^{(n)} - \theta \bigr|
	\Bigl|l_k(\theta-\tfrac{\pi}{2n}) +
	l_k(\theta+\tfrac{\pi}{2n}) \Bigr|\delta^{-1}
	\,\omega(f,\delta).
\end{aligned}
\]
Let 
\[
\delta \;=\delta_n=
\frac{1}{2}\sum_{k=1}^n 
\bigl|\theta_k^{(n)} - \theta \bigr|
\Bigl|l_k(\theta-\tfrac{\pi}{2n}) +
l_k(\theta+\tfrac{\pi}{2n}) \Bigr|\leq \left(\frac{\pi}{2n}+\epsilon_n\right)c_2
+
\frac{9\pi^4}{64}\frac{1}{n\epsilon_n^2}.
\]
The last inequality has been derived in Theorem \ref{voronov}, where $\epsilon_n=\frac{1}{n^{\frac{1}{3}}}$. Let 
\[
\mu_n:=\left(\frac{\pi}{2n}+\epsilon_n\right)c_2
+
\frac{9\pi^4}{64}\frac{1}{n\epsilon_n^2}=\left(\frac{\pi}{2n}+\frac{1}{n^{\frac{1}{3}}}\right)c_2
+
\frac{9\pi^4}{64}\frac{1}{n^{\frac{1}{3}}}.
\]  
Hence, $|\mathcal{G}_n(f)(\theta)-f(\theta)|\leq (c_2+1)\omega(f,\mu_n)$ for all $\theta\in [0,\pi]$. The result follows.

\end{proof}

	\section*{Concluding Remarks and future problems}
	In this article, we defined a sequence of operators \( \{\mathcal{G}_n\}\) using the classical Lagrange interpolation operators. Although it is well known that the Lagrange interpolation operators with respect to the Chebyshev nodes does not converge, by a convergence result established by G. Grünwald in $1941$, the sequence \( \{\mathcal{G}_n(f)\}\) converges uniformly to \(f\) for all $f$ in  $C[0,\pi]$, with respect to the Chebyshev nodes. Following this, we established a perturbed version of the result, thereby extending the class of nodal sets for which uniform convergence holds. Then we discussed some cases where similar convergence holds. \par
	The operators $\mathcal{G}_n$ are non-positive for $n=1,2,\ldots$ . On the Chebyshev nodes, this has been obtained in \cite{vin}. This technique can be adapted for general nodes.\par 
	It is known that the Lagrange interpolation operator based on equidistant nodes diverges; see \cite{mohapatra, revers}. This naturally raises the question: Does the averaged operator \( \mathcal{G}_n \), defined via a symmetric perturbation of the interpolation points converge when the nodes are equidistant? Can we characterize the class of nodal sets \( \{ \eta_k \}_{k=1}^n \) for which the sequence \( \mathcal{G}_n(f) \) converges uniformly for every \( f \in C[0,\pi] \)?
	
	\section*{Acknowledgment}
	The author sincerely thanks Prof.~M.~N.~N.~Namboodiri for suggesting this interpolation problem. She also expresses her gratitude to Dr.~V.~B.~Kiran Kumar and Prof.~M.~N.~N.~Namboodiri for fruitful discussions and valuable suggestions which contributed to the development of this work.
	\section*{Funding} 
	No funds were received for carrying out this research.
	\bibliographystyle{elsarticle-harv}

\end{document}